# Social Choice Under Incomplete, Cyclic Preferences

## Majority/Minority-Based Rules, and Composition-Consistency

**Jobst Heitzig**

Institut für Mathematik, Universität Hannover, Welfengarten 1, D-30167 Hannover, Germany, e-mail: `heitzig@math.uni-hannover.de`, fax: +49 511 762 5803

**Abstract** Actual individual preferences are neither complete (=total) nor anti-symmetric in general, so that at least every *quasi-order* must be an admissible input to a satisfactory choice rule. It is argued that the traditional notion of "indifference" in individual preferences is misleading and should be replaced by *equivalence* and *undecidedness*.

In this context, ten types of *majority* and *minority arguments of different strength* are studied which lead to social choice rules $\mathcal{C} : (S; \mathbb{R}) \mapsto C(S; \mathbb{R}) \in \mathcal{P}(S) \setminus \{\emptyset\}$ that accept profiles $\mathbb{R}$ of *arbitrary* reflexive relations. These rules are discussed by means of many familiar, and some new conditions, including *immunity from binary arguments*. Moreover, it is proved that every choice function satisfying two weak Condorcet-type conditions can be made both composition-consistent and idempotent, and that all the proposed rules have polynomial time complexity.

## 1 Introduction

In this paper, we shall reconsider the question of how a group of individuals should proceed to make their decision if they want to realize exactly one out of a certain set $X$ of alternatives. We shall do this under the usual assumptions that (i) all information relevant for the choice is which individuals prefer which alternatives to which other alternatives, and that (ii) the result of the procedure should not only be some alternative(s) that are "best" in $X$, but rather some "best" alternative(s) for each non-empty subset $S$ of $X$. For the sake of simplicity, we will not consider "degrees" of preference here, although there are some good arguments that, for

I would like to thank the editors for calling my attention to articles [**?**] and [**?**].

example, *fuzzy* relations are a more precise model of actual individual preferences. Still, the ideas presented here can probably be translated into a fuzzy preference setting, particularly, the classes of generalized majority/minority relations we will study below may perhaps be interpreted as fuzzy *social* preferences. But although the introduced axioms of *immunity from binary arguments* can be interpreted as axioms for one individual with fuzzy preferences (as in [**?**], for example) rather than a society with exact preferences, they are not intended for and hardly make sense in the fuzzy case.

If we understand (i) in a way that denies either the existence of "degrees" of individual preference, or their measurability, or at least their relevance, we can assume that the individual preferences are given as binary relations on $X$, and since their significance lies in the comparison of *different* alternatives, we may adopt the convention to use reflexive relations only. We will see, however, that it is neither justified nor (for most purposes) necessary to assume that individual preference relations have any special properties other than reflexivity.

As usual, we deal with a finite set $N$ of individuals, and without loss of generality let $N = \{1, \dots, n\}$ and $n \geqslant 1$. I will use $i, j, \dots$ as variables over individuals. The possible alternatives build a set $X$, and it is important that $X$ is finite, too, but this is not an actual restriction since, in all practical situations, there will always be only finitely many feasible alternatives. We may as well assume that $|X| = m \geqslant 2$. Alternatives will be denoted by variables $x, y, \dots$.

The motivation for (ii) is this: Because $X$ is meant to contain also those alternatives that might turn out to be actually impossible after the decision is made, and because, on the other hand, there are situations in which two or more alternatives appear in the individual preferences in completely equivalent ways, we shall consider *social choice rules* $\mathcal{C}$, i.e. algorithms that provide the group with a (multivalued) *choice function* rather than a single alternative. This is a function $C : \mathcal{P}(X) \setminus \{\emptyset\} \to \mathcal{P}(X) \setminus \{\emptyset\}$ that assigns to each nonempty $S \subseteq X$ a (hopefully small) subset $C(S) \subseteq S$ of alternatives that I will call *acceptable* here. Then, as soon as the set $S$ of *feasible* alternatives is known to the group, they can choose some $x \in C(S)$ either randomly, or by delegating that decision to some individual, or in whatever alternative way. In fact, some authors discuss algorithms that directly assign a probability distribution to $S$ rather than a subset. But even in such a model, which is more general than the present one, a discussion of choice rules remains important, since it is rather natural to define an *induced* choice rule by letting $C(S)$ be the set of alternatives receiving a *nonzero* probability.

In order to have a fixed interpretation of the input to the algorithm $\mathcal{C}$, let us assume that, for each pair $(x, y) \in X \times X$ each individual $i \in N$ has been

> *asked whether or not she thinks that*
> *alternative $x$ is at least as desirable to her as alternative $y$.*

The input to the algorithm is then the resulting *profile (of individual preferences) on $X$*, i.e. the tuple $\mathbb{R} = (R_1, \dots, R_n)$ of reflexive binary relations on $X$, where $x\, R_i\, y$ holds if and only if $i$ has answered "yes" to the above question. To make clear from which profile $\mathbb{R}$ the output choice function is derived, the set $C(S)$ will in some places more accurately be denoted by $C(S; \mathbb{R})$.

An *isomorphism* between profiles $\mathbb{R}$ on $X$ and $\mathbb{R}'$ on $X'$ is here a bijection $\varphi : X \to X'$ between the two sets of alternatives for which there is another bijection $\psi : N \to N'$ between the two sets of individuals such that $x\,R_i\,y \iff \varphi(x)\,R'_{\psi(i)}\,\varphi(y)$ for all $x, y \in X$ and $i \in N$. In this case, $\mathbb{R}$ and $\mathbb{R}'$ will be called *isomorphic*. For example, the identity map $id_X$ is an isomorphism between $\mathbb{R} = (R_1, \ldots, R_n)$ and $\mathbb{R}^\psi := (R_{\psi(1)}, \ldots, R_{\psi(n)})$ for every permutation $\psi : N \to N$ of individuals. Throughout this paper I will adopt a very broad idea of equality and independence: We will only consider rules that are *anonymous* w.r.t. individuals and *neutral* w.r.t. alternatives, which can be summarized in the following condition of *isomorphism invariance:*

(Iso) $C(\varphi[S]; \mathbb{R}') = \varphi[C(S; \mathbb{R})]$
whenever $\varphi$ is an isomorphism between $\mathbb{R}$ and $\mathbb{R}'$.

We may then think of $\mathbb{R}$ and $\mathbb{R}^\psi$ as being essentially the same profile. Moreover, our rules will be *independent of irrelevant alternatives,* i.e. fulfill

(I) $C(S; \mathbb{R}) = C(S; \mathbb{R}|_S)$,

where $\mathbb{R}|_S = (R_1|_S, \ldots, R_n|_S)$ and $R_i|_S = R_i \cap (S \times S)$ are the *restrictions* of $\mathbb{R}$ and $R_i$ to $S$, respectively. The condition expresses the idea that, for a choice of feasible alternatives, only (preferences about) feasible alternatives should be relevant. This has the nice consequence that $\mathcal{C}$ is already determined as soon as $C(X; \mathbb{R})$ is known for all $X$ and all profiles $\mathbb{R}$ on $X$.

## 2 At least all quasi-orders must be admissible individual preferences

To begin with, a given profile $\mathbb{R}$ provides us with the following additional relations: $P_i := R_i \setminus R_i^{op}$ is the asymmetric part of $R_i$, where $R_i^{op} = \{(x, y) : (y, x) \in R_i\}$, and *asymmetric* means that $x\,P_i\,y$ excludes $y\,P_i\,x$. $P_i$ contains the *expressed strict preferences* of $i$. The rest of $R_i$ is its symmetric part $E_i := R_i \cap R_i^{op}$, the *expressed equivalences* of $i$. Although, in case of complete preferences, this relation is usually called "indifference" and is then consequently denoted by the letter $I$, in case of incomplete preferences the above terminology is better for two reasons: It resembles standard order-theoretical terminology (cf. [?]), and, what is especially important in this context, $E_i$ must not be confused with the relation $x\,U_i\,y :\iff \neg(x\,R_i\,y \vee y\,R_i\,x)$ which, according to the interpretation we started with, encodes the *expressed undecidedness* of $i$ whether to prefer $x$, or $y$, or whether to value them as equally desirable.

Being undecided does not automatically imply being also *unconcerned*, because it may simply be the case that a person measures the alternatives with more than just one criterion and that $x$ is better with respect to one aspect but $y$ with respect to another. My point is that, in fact, in the fewest situations each individual uses only one criterion to evaluate alternatives. Let us therefore consider the situation of an individual $i$ with a set $C_i$ of different criteria so that each criterion $c \in C_i$ provides him with a different preference relation $R_{ic}$. There is no legitimation (and no necessity) to forbid him to consider $x$ at least as desirable as $y$ if and only if $x$ is at least as good as $y$ with respect to *all* his criteria. More precisely: It well makes sense and should therefore surely be legitimate for $i$ to represent his preferences by the relation $R_i := \bigcap_{c\in C_i} R_{ic}$. It should also be his right to use as many and as independent criteria as he likes. But then the resulting $R_i$ might be any *partial order* (i.e. reflexive, transitive, antisymmetric, but *not* necessarily complete), even if all used criteria lead to *linearly ordered* (i.e. complete,[1] transitive, and antisymmetric) preferences $R_{ic}$. This follows from a basic order-theoretical fact (cf. [**?**]):

**Lemma 1** *Any partial order on a finite set $X$ is an intersection of at most $|X|$ many linear orders on $X$, any quasi-order on $X$ an intersection of at most $|X|$ many total quasi-orders of the form $(X \times X) \setminus \big((X \setminus S) \times S\big)$ for some $S \subseteq X$.*

In the latter kind of total quasi-orders there are two subsets $S$, $X \setminus S$ of equally desirable alternatives such that those from $S$ are preferred to the rest. They naturally arise from *binary* criteria, therefore Lemma 1 implies that even if only binary criteria are used, $R_i$ might still be even any *quasi-order* (i.e. reflexive and transitive, but not necessarily complete or antisymmetric). The immediate and in my eyes indisputable consequence of this is that, contrary to the completeness assumption that occurs throughout the literature, *every quasi-order on $X$ must be an admissible individual preference relation.*

A quite different argument against the completeness assumption comes from the fact that $i$ might have different intentions when expressing $x\,E_i\,y$ rather than $x\,U_i\,y$: On the one hand, saying "$x$ and $y$ are equally desirable to me" can be interpreted as a vote not to distinguish between $x$ and $y$, i.e., a vote for having either both or none of them in $C(S)$. On the other hand, $x\,U_i\,y$ may be interpreted as the statement "I do not *want* to decide about $x$ and $y$" in the sense that $i$ wants to *delegate* the decision about $x$ and $y$ to those individuals that have more information about, or more interest in the distinction between $x$ and $y$. Such a delegation is an often sensible and, in practice, very common behaviour, especially when (i) there are many alternatives, or (ii) some pairs of alternatives differ only in their effect on few individuals, or (iii) some individuals have restricted information.

Before we turn to the concept of majority and minority, let us finally remember that it has also been noticed that, in some cases, individual preferences may contain cycles. But, despite the fact that we have to be a bit careful with Pareto-type principles then (see below), there seems to be no problem at all in dealing with cyclic preferences, too, at least not when the preference information that is taken into account by the rule $\mathcal{C}$ in any case only consists of some cardinalities like the following:

$$\begin{aligned}
r_{xy} &:= |\{i \in N : x\,R_i\,y\}|, \\
p_{xy} &:= |\{i \in N : x\,P_i\,y\}|, \\
e_{xy} &:= |\{i \in N : x\,E_i\,y\}| = e_{yx} = r_{xy} - p_{xy}, \\
u_{xy} &:= |\{i \in N : x\,U_i\,y\}| = u_{yx} = n - d_{xy},
\end{aligned}$$

---

[1] A note on terminology: in order theory, complete relations are usually called "total" instead, while a "complete" quasi-order is one in which infima and suprema exist.

$$d_{xy} := |\{i \in N : x\,R_i\,y \lor y\,R_i\,x\}| = d_{yx} = p_{xy} + p_{yx} + e_{xy}, \text{ and}$$
$$d_S := |\{i \in N : x\,P_i\,y \text{ for some } x, y \in S\}|.$$

In some places below, I will identify irreflexive relations with their reflexive counterparts and use a somewhat sloppy arrow notation; for example, $A = a \rightleftarrows c \to b \leftarrow a, d \to b$ should be read as

$$A = \{(a,a), (b,b), (c,c), (d,d), (a,c), (c,a), (c,b), (a,b), (d,b)\}.$$

Only in case of quasi-orders, the more convenient Hasse diagrams can, and will be used; the above example would therefore rather look like this:

$$A = \begin{matrix} ac & d \\ \searrow & \swarrow \\ & b \end{matrix}.$$

Analogously to the order-theoretic notions of "minimal" and "smallest" elements, let us call an alternative $x$ *A-optimal* if, for all $y$, $y\,A\,x$ implies $x\,A\,y$. As there are not always optimal elements of $A|_S$, we will often use its *transitive hull* $\mathrm{tr}_S(A) = \bigcup_{k=0}^{\infty}(A|_S)^k$, that is, the smallest quasi-order on $S$ containing $A|_S$. If $A$ is a *tournament*, that is, complete and antisymmetric, the set of $\mathrm{tr}_S(A)$-optimal elements is just the top-cycle of $A|_S$.

## 3 Binary arguments supported by majorities

Despite some irritating phenomena that are related to the concept of "majority", this notion is surely the most important idea in the theory of social choice. It is well known that we can't expect any alternative to have majority support, but whatever exact definition of majority we may adopt, an alternative $x$ should only be acceptable to the group if it is at least in some sense "defendable" against arguments a majority of the group might give *against* $x$. The most important type of such arguments seems to be what may be called "binary" arguments: A part of the group might argue that another alternative $y$ is more desirable to them than $x$ and that therefore $x$ should not be acceptable. Since, for every alternative there might be a majority that favours a different alternative, we should think about possibilities to counter and *refuse* some of these binary arguments so as to make a choice possible. Assume that we have in some way decided which kinds of binary arguments to consider important, and let $y\,A\,x$ denote the fact that the given profile $\mathbb{R}$ contains such an argument for $y$ against $x$. Now consider the following possible conditions on $\mathcal{C}$:

(wIm$_A$) *If $x \in C(S)$, $y \in S \setminus \{x\}$, and $y\,A\,x$,*
  *there must be $z \in S \setminus \{y\}$ with $z\,A\,y$.*

(Im$_A$) *If $x \in C(S)$, $y \in S$, and $y\,A\,x$, then $x\,\mathrm{tr}_S(A)\,y$*
  *(i.e. there must be $z_1, \ldots, z_m \in S$ with $x\,A\,z_1\,A\cdots A\,z_m\,A\,y$).*

(sIm$_A$) *Each $x \in C(S)$ must be optimal in $\mathrm{tr}_S(A)$*
  *(i.e. if $x \in C(S)$, $y \in S$, and $y\,\mathrm{tr}_S(A)\,x$, then also $x\,\mathrm{tr}_S(A)\,y$).*

The condition (wIm$_A$) of *weak A-immunity* claims that when $y$ is used in an argument against $x$, $x$ can only be acceptable if the proposed "better" alternative $y$ is subject to a similar argument.

Perhaps one should at least require the stronger property of *A-immunity* (Im$_A$) which claims that the argument $y\,A\,x$ must even be answered with a sequence of similar arguments that lead back to $x$. This is because the existence of such a sequence effectively demonstrates that the argument $y\,A\,x$ is destructive in two respects: (i) It cannot "consistently" be taken into account without the risk of making *all* alternatives inacceptable (rather than just $x$): if the argument $y\,A\,x$ is successful, then so should be all others in the sequence, which would result in excluding $x, y, z_1, \ldots, z_m$ at once. (ii) The argument is also somehow useless for its supporters, because it does not place $y$ in an essentially better position than $x$.

The appeal of the even more restrictive *strong A-immunity* (sIm$_A$) is that it treats the situation more symmetrically: Even if there is only a sequence of arguments leading from $y$ to $x$ instead of a direct one, there shall also be such a sequence leading back.

Let us look at some possible concretizations of the notion of majority, that is, some possible definitions of $A$ as a function of $\mathbb{R}$. The following binary relations on $X$ encode (proper) majorities of different types and strengths $\alpha$ in terms of the individual preferences: For all $S \subseteq X$, $x, y \in S$, and $\frac{1}{2} < \alpha \leqslant 1$, let

$$x\,M_\alpha\,y :\iff p_{xy} \geqslant \alpha n, \qquad (\implies p_{xy} > 0)$$
$$x\,N_\alpha\,y :\iff r_{xy} \geqslant \alpha n, \qquad (\implies r_{xy} > 0)$$
$$x\,M_\alpha^S\,y :\iff p_{xy} \geqslant \alpha d_S > 0,$$
$$x\,N_\alpha^S\,y :\iff r_{xy} \geqslant \alpha d_S > 0,$$

$$x\,B_\alpha\,y :\iff p_{xy} \geqslant \alpha d_{xy} > 0,$$
$$x\,D_\alpha\,y :\iff r_{xy} \geqslant \alpha d_{xy} > 0,$$

$$x\,P_\alpha\,y :\iff p_{xy} \geqslant \alpha(p_{xy} + p_{yx}) > 0,$$
$$x\,R_\alpha\,y :\iff r_{xy} \geqslant \alpha(r_{xy} + r_{yx}),\ p_{xy} > 0$$

$$x\,U_\alpha\,y :\iff p_{xy} > p_{yx},\ p_{xy} + u_{xy} \geqslant \alpha n$$
$$\qquad\qquad (\iff r_{xy} > r_{yx},\ r_{yx} \leqslant (1-\alpha)n),\text{ and}$$
$$x\,E_\alpha\,y :\iff p_{xy} > p_{yx},\ p_{xy} + e_{xy} \geqslant \alpha d_{xy}$$
$$\qquad\qquad (\iff r_{xy} > r_{yx},\ r_{xy} \geqslant \alpha d_{xy}).$$

(For $\alpha = 1$, $P_\alpha$, $R_\alpha$, $U_\alpha$ and $E_\alpha$ should, of course, not be confused with the preferences of individual 1.)

Most of these definitions have in common that the individuals constituting the majority share some opinion about $x$ and $y$ and build a fraction of at least $\alpha$ of all in some way relevant individuals. Moreover, there must at least be one individual within a majority. That unifying opinion can be either *strictly preferring $x$ to $y$* (which corresponds to using the number $p_{xy}$ in $M_\alpha$, $B_\alpha$, and $P_\alpha$, for example), or just *considering $x$ at least as desirable as $y$* (which corresponds to using $r_{xy}$ in $N_\alpha$, $D_\alpha$, and $R_\alpha$). The relevant individuals are either all (where $n$ is used), or those that are not undecided about $x$ and $y$ (which corresponds to using $d_{xy}$), or only those that express a strict preference about $x$ and $y$ (as in $P_\alpha$). Although $r_{xy} + r_{yx}$ does not, in general, enumerate some subset of individuals, $R_\alpha$ has the appeal that, similarly to $P_\alpha$, it considers those majority arguments "equally strong" in which the proportion $r_{xy}/r_{yx}$ is constant. This is because $x\,R_\alpha\,y$ is mainly equivalent to $r_{xy}/r_{yx} \geqslant \frac{\alpha}{1-\alpha}$ (the condition $p_{xy} \neq 0$ is added only to ensure antisymmetry and the inclusions that are stated below). In some of the definitions, $d_S$ occurs instead of $n$, which models kinds of "semi-relative" majorities. This could provide a compromise between using absolute majorities and the possible requirement that $C(S)$ should be at least independent of individuals that are completely undecided about all feasible alternatives. On the other hand, it might be problematic if whether or not individual $i$ is counted would depend on perhaps just one alternative's being feasible or not.

The last two definitions deserve a more detailed explanation: They express the idea that in the first place only individuals with strict preferences for $x$ over $y$ will *raise* the corresponding argument against $y$, but that then additional individuals may *join* them to build a majority. For example, if only slightly more people strictly prefer $x$ to $y$ than $y$ to $x$, the former might persuade all who *equally* desire $x$ and $y$ to support their argument, so that in the end a very strong majority evolves (relative to all decided individuals); this is assumed in the definition of $E_\alpha$. One could also assume that, on the contrary, the *undecided* individuals are persuaded to constitute an (absolute) majority, which leads to $U_\alpha$.

Since $\alpha > \frac{1}{2}$, it can easily be seen that, despite $N_\alpha$, $N_\alpha^S$, and $D_\alpha$, all of the above relations are antisymmetric. Moreover, we have some, perhaps unexpected, inclusions that are shown in the lower part of the diagram on page 10: (i) When $\alpha$ grows, the relations obviously shrink (which is represented by the dotted lines).

$$\text{(ii)}\ M_\alpha \subseteq M_\alpha^S \subseteq B_\alpha \subseteq R_\alpha \subseteq P_\alpha \subseteq E_\alpha \subseteq D_\alpha,$$

because $n \geqslant d_S \geqslant d_{xy}$, $r_{xy} = p_{xy} + e_{xy}$, and

$$\alpha d_{xy} + e_{xy} \geqslant \alpha(r_{xy} + r_{yx}) \geqslant \alpha(p_{xy} + p_{yx}) + e_{xy} \geqslant \alpha d_{xy}.$$

Similarly, $M_\alpha \subseteq N_\alpha \subseteq N_\alpha^S \subseteq D_\alpha$, $M_\alpha^S \subseteq N_\alpha^S$, and $B_\alpha \subseteq U_\alpha$. (iii) For $\alpha \leqslant \frac{n}{2n-1}$, $x\,R_\alpha\,y$ is equivalent to $p_{xy} \neq 0$ and $r_{xy}/r_{yx} \geqslant \frac{n}{n-1}$, the latter being equivalent to $r_{xy} > r_{yx}$. Then $x\,R_\alpha\,y$ is already implied by $x\,E_\alpha\,y$, so that $R_\alpha = P_\alpha = U_\alpha$. (iv) $x\,U_\alpha\,y$ implies $\neg\,y\,N_{1-\alpha+1/n}\,x$.

Summarizing this section, we have the following classes of majority relations: $M_\alpha$ and $N_\alpha$ encode *absolute* majorities, $M_\alpha^S$ and $N_\alpha^S$ *semi-relative* majorities, $B_\alpha$ and $D_\alpha$ *relative* ones, $P_\alpha$ and $R_\alpha$ *proportional* ones, and $E_\alpha$ and $U_\alpha$ encode *persuaded* majorities. Moreover, $M$, $M^S$, $B$, $P$, $E$ and $U$ will be called the *strict* types, and the rest *non-strict*.

## 4 Rules based on immunity from classes of binary arguments

For a fixed profile $\mathbb{R}$, each of the sets $\{\text{tr}_S(A_\alpha) : \frac{1}{2} < \alpha \leqslant 1\}$, where $A$ is one of the ten types $M$, $N$, $M^S$, ..., is a chain of quasi-orders. Now it is important that $X$ is finite: Then the chain is also finite, and the following Lemma applies:

**Lemma 2** *Every chain[2] $Q_1 \subseteq \cdots \subseteq Q_m$ of finite quasi-orders has a common optimal element, i.e. some $x$ that is $Q_k$-optimal for all $k \leqslant m$.*

*Proof.* Since $Q_m$ is finite, the set $S_m$ of its optimal elements is not empty. Let $S_{m-1} \subseteq S_m$ be the (also non-empty) set of optimal elements of $Q_{m-1}|_{S_m}$. Then each $x \in S_{m-1}$ is also $Q_{m-1}$-optimal, because

$$y\, Q_{m-1}\, x \Longrightarrow y\, Q_m\, x \Longrightarrow y \in S_m \Longrightarrow x\, Q_{m-1}\, y.$$

Thus, denoting the set of optimal elements of $Q_k|_{S_{k+1}}$ by $S_k$, we inductively get $S_1$, the still not empty set of all common optimal elements. □

Note that, consequently, the set of common optimal elements can be found in polynomial time. The above observation enables us to fulfill (sIm$_A$) not only for one specific majority relation, say $M_{1/2+\varepsilon}$, but at once for a whole class of majority relations of the same type but of different strength: We may simply define $C(S)$ as the set of common optimal elements of, for example, the chain $\{\text{tr}_S(M_\alpha) : \frac{1}{2} < \alpha \leqslant 1\}$. Then any $x \in C(S)$ can be defended against a sequence of arguments $y\, M_\alpha\, z_1\, M_\alpha \cdots M_\alpha z_m\, M_\alpha\, x$ with a sequence of *equally strong* $M$-type arguments. Exchanging $M$ by any of the other types of majority relations and varying the lower and upper bounds for $\alpha$, we get a large number of *majority-based rules* with very good immunity properties. Let us introduce the following notation: For $\frac{1}{2} \leqslant \beta < \gamma \leqslant 1$, the rule based on the chain $\{\text{tr}_S(A_\alpha) : \beta < \alpha \leqslant \gamma\}$, where $A$ is one of $M$, $N$, $M^S$, ..., will be denoted by $\mathcal{A}_{(\beta,\gamma]}$. For example, we just introduced the rule $\mathcal{M}_{(.5,1]}$.

It seems quite natural to me to take *minorities* in the same way into account as majorities: Suppose we apply the definitions of $M_\alpha, \ldots, U_\alpha$ even for $0 < \alpha \leqslant \frac{1}{2}$, except that we drop all of the " $> 0$"-requirements here. Then $D_\alpha$, $P_\alpha$, and $R_\alpha$ become complete relations for $\alpha \leqslant \frac{1}{2}$ (but only $U_\alpha$ and $E_\alpha$ remain antisymmetric in general), and in the inclusions, $R_\alpha$ and $P_\alpha$ change places: Still $M_\alpha \subseteq N_\alpha \subseteq N_\alpha^S \subseteq D_\alpha$, but now

$$M_\alpha \subseteq M_\alpha^S \subseteq B_\alpha \subseteq P_\alpha \subseteq R_\alpha \subseteq D_\alpha$$

and $U_\alpha \subseteq E_\alpha = R_{n/(2n-1)}$. Moreover, $R_{1/2} = P_{1/2}$, and no longer $B_\alpha \subseteq U_\alpha$.

Now, requiring (Im$_{A_\alpha}$) as well for $0 < \alpha \leqslant \frac{1}{2}$ would ensure that also a minority argument $y\, A_\alpha\, x$ will be successful if only if it can be taken into account "consistently". In such a case it is not only harmless to give the minority this "direct" power to exclude $y$, but it may even be *indicated* in certain situations: (i)

---

[2] Chains (towers/nested families) of other kinds of relations are also an important tool in modelling uncertain or stochastically varying preferences of one individual, cf. [**?**,**?**]

When an absolute or semi-relative majority type is used and many individuals are undecided (because of too few information, for example), or (ii) when a strict type other than $P$ is used and many individuals equally desire $x$ and $y$ (because they are not affected by their distinction, for example), minority arguments should be taken into account.

Moreover, if we even require $(\text{sIm}_{A_\alpha})$ for $0 < \alpha \leqslant \frac{1}{2}$, a minority argument $y\,A_\alpha\,x$ which is refused as a *single* argument, may still have "indirect" power: Other such arguments can help constituting a *sequence* of arguments leading to the exclusion of $x$, if at least one of them is non-refusable. While indirect influence works with all of the types, direct influence of a minority is only possible if $A$ is none of the complete types $D$, $P$, or $R$: If, for example, $y\,D_\alpha\,x$ for some $\alpha \leqslant \frac{1}{2}$ but for no $\alpha' > \frac{1}{2}$, then $r_{yx} \leqslant \frac{1}{2}d_{yx}$, thus $r_{xy} \geqslant \frac{1}{2}d_{xy} \geqslant \alpha d_{xy}$, hence also $x\,D_\alpha\,y$. However, in the crucial situations (i) and (ii) described above, a subgroup of individuals that would be a minority of one of the other types might well constitute a majority of type $D$, $P$, or $R$.

The rules $\mathcal{A}_{(\beta,\gamma]}$ with $0 \leqslant \beta \leqslant \gamma \leqslant \frac{1}{2}$ will be called *minority-based rules*, those with $0 \leqslant \beta \leqslant \frac{1}{2} < \gamma \leqslant 1$ *mixed rules*.

As for the question which of the possible types of majority/minority to use, it may turn out that this cannot be completely decided in the usual axiomatic way. Thus, in addition to the short axiomatic discussion in the next section, the rules should also be compared from a more practical perspective. Surely, some types like $M$, $N$, $B$, $P$, and $D$ look more "simple" or "natural" than others, and $E$ and $U$ are based on questionable behavioural imputations. There is a good reason to prefer some of the "larger" types: Whatever type $A$ we use, an $x \in C(S)$ may become subject to an argument supported by a majority of a different type, say $y\,A'_\alpha\,x$. It may then still possible to refuse $y\,A_\alpha\,x$ on the basis of a sequence of $A_\alpha$-arguments leading back to $y$, if only $A'_\alpha \subseteq A_\alpha$, i.e., the rule using $A$ is in some sense also immune from arguments corresponding to types smaller than $A$.

Following an argument from Section 2, another criterion is that $x\,E_i\,y$ and $x\,U_i\,y$ should not always have the same effect, which is a point against using $M$ and $P$.

In all, the non-persuaded, relative types $D$ (being the largest such) and $B$ (being the largest that allows "direct" influence of minorities) seem to give the best compromises so far. As for the range of $\alpha$, its lower bound $\beta$ should be taken as small as possible in order to keep the resulting set $C(S)$ as small as possible. On the other hand, taking the upper bound $\gamma$ not to large could be a sort of protection of minorities, since then it would be possible to counter arguments of strength $> \gamma$ by a sequence of arguments of strength $\gamma$. However, even in such a case one would probably add $\text{tr}_s(A_1)$ to the chain, the corresponding rules will be denoted by $\mathcal{A}_{(\beta,\gamma],1}$ Depending on how much power and protection we want to give minorities, one could for example use the rules $\mathcal{D}_{(.5,1]}$, $\mathcal{D}_{(0,1]}$, $\mathcal{B}_{(0,1]}$, $\mathcal{B}_{(0,2/3],1}$, or even $\mathcal{B}_{(0,.5],1}$.

**Fig. 1** Majority and minority relations

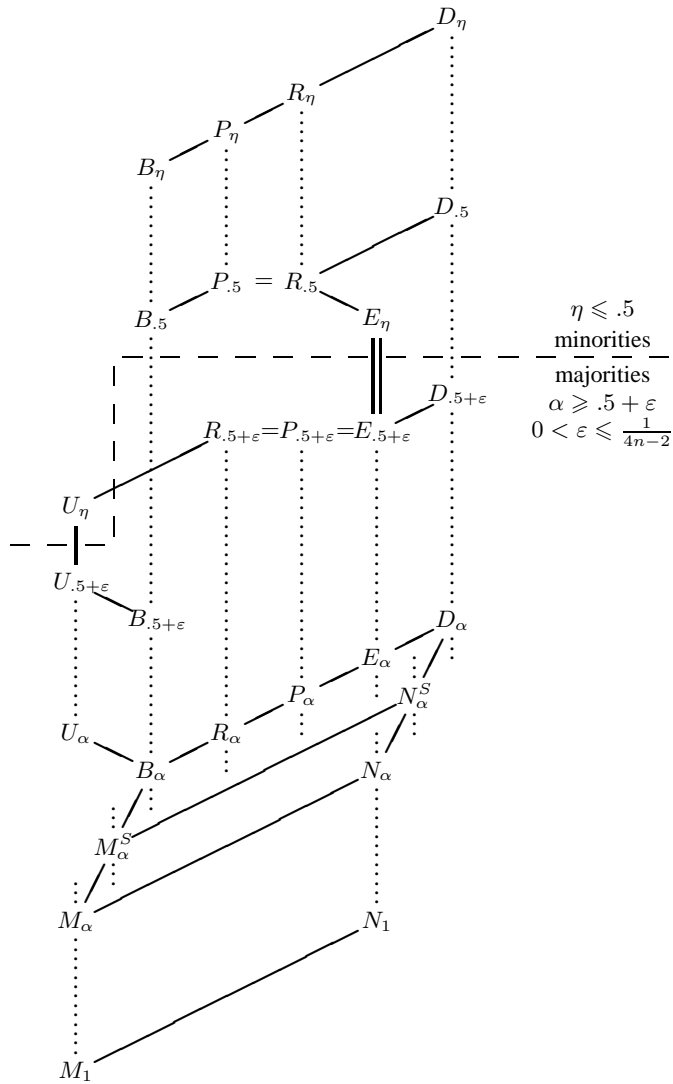

## 5 Some conditions discussed

**Pareto-principles.** The idea that $y$ is unacceptable if all individuals (strictly) prefer $x$ to it can only work when there is a minimal amount of rationality in at least one individual's preferences. Thus, in a setting where cyclic individual preferences are explicitly admitted, adequate formulations of this idea must take possible cycles into account:

(wP) $y \notin C(S)$ whenever $x\,M_1\,y$ for some $x \in S$, but not $y \operatorname{tr}_S(M_1)\,x$.

(sP) $y \notin C(S)$ whenever $x\,N_1\,y$ for some $x \in S$, but not $y \operatorname{tr}_S(N_1)\,x$.

(sP') $y \notin C(S)$ whenever $x\,B_1\,y$ for some $x \in S$, but not $y \operatorname{tr}_S(B_1)\,x$.

(NNP) $x \in C(S)$ whenever $x \in S$ and $x\,N_1\,y$ for some $y \in C(S)$.

The *weak* and *strong Pareto*-principles (wP) and (sP) are just the conditions $(\operatorname{Im}_{M_1})$ and $(\operatorname{Im}_{N_1})$, while the alternative (sP') is just $(\operatorname{Im}_{B_1})$. Since $M$ is the smallest majority type, adding the relation $\operatorname{tr}_S(M_1)$ to a chain $\{\operatorname{tr}_S(A_\alpha)\}$ gives again a chain. This results in a class of (wP)-*modified* rules I will denote by $\mathcal{A}_{(\beta,\gamma]}\mathcal{M}_1$. The case of (sP) is different: $\operatorname{tr}_S(N_1)$ cannot always be added, but:

**Lemma 3** *Suppose $y$ is $\operatorname{tr}_S(A_\alpha)$-optimal, and $x\,N_1\,y$ for some $x \in S$. Then $x$ is also $\operatorname{tr}_S(A_\alpha)$-optimal in each of the following cases:*

*(i) $A \in \{N, N^S, D\}$.*
*(ii) $A \in \{M, M^S, P\}$, and $R_i P_i R_i \subseteq P_i$ for all $i \in N$.*
*(iii) $A \in \{B, R, E\}$, and all $R_i$ are transitive.*

*Proof.* Immediate since (i) $N_1 \subseteq N_\alpha \cap N_\alpha^S \cap D_\alpha$, (ii) $x\,N_1\,y$ and $R_i P_i R_i \subseteq P_i$ for all $i \in N$ imply $p_{xz} \geqslant p_{yz}$ and $p_{zy} \geqslant p_{zx}$ for all $z \in X$, and (iii) $x\,N_1\,y$ and $R_i R_i \subseteq R_i$ for all $i \in N$ imply $r_{xz} \geqslant r_{yz}$ and $r_{zy} \geqslant r_{zx}$. $\square$

Therefore, among the common $\operatorname{tr}_S(A_\alpha)$-optimal elements there is at least one $N_1$-optimal element as long as (i), (ii), or (iii) holds. Particularly, we get the (sP)-*modified* rules $\mathcal{N}_{(\beta,\gamma],1}$, $\mathcal{N}_{(\beta,\gamma]}^S\mathcal{N}_1$, and $\mathcal{D}_{(\beta,\gamma]}\mathcal{N}_1$ by adding $\operatorname{tr}_S(N_1)$ to the defining chain. For other types $A$, the analogously defined algorithms $\mathcal{A}_{(\beta,\gamma]}\mathcal{N}_1$ only produce a nonempty set $C(S)$ if all individual preferences show at least some minimal amount of "rationality" in form of some transitivity property. On the other hand, in my opinion, (sP) is not always reasonable anyway: It misinterprets equivalence as undecidedness and ignores that $x\,E_i\,y$ often expresses $i$'s intention that either both or none of $x$ and $y$ be accepted. This is of particular importance if the alternatives at hand affect individuals that do not belong to the deciding group, as for example when an award must be given to a single candidate: Suppose that, except one, all members $i$ of a not too small jury have $x\,E_i\,y$. Then it seems wise to rather select one of the candidates $x, y$ randomly than excluding one of them only on the basis of a single juror's preferences.

Analogously, the second strong form (sP') of the Pareto-principle can be fulfilled by using an (sP')-*modified* rule $\mathcal{A}_{(\beta,\gamma]}\mathcal{B}_1$ with $A \in \{B, U, R, P, E, D\}$. Again by the previous lemma, $\mathcal{A}_{(\beta,\gamma]}$ fulfills the last condition (NNP), *non-negative Pareto*, for all profiles if $A \in \{N, N^S, D\}$, but only for "rational" profiles if $A \in \{M, M^S, P, B, R, E\}$. Finally, for profiles of quasi-ordered preferences, one may satisfy all four conditions using an algorithm $\mathcal{A}_{(\beta,\gamma]}\mathcal{B}_1\mathcal{M}_1\mathcal{N}_1$ with $A \in \{B, R, P, E, D\}$.

**Responsiveness.**

(sNNR) If $x \in C(S;\mathbb{R}) \not\ni z$, $y \in X$, $i \in N$, $R_i' \setminus \{x,y\}^2 = R_i \setminus \{x,y\}^2$,
and either $x\,E_i\,y$ and $x\,P_i'\,y$, or $y\,P_i\,x$ and $y\,E_i'\,x$,
then $x \in C(S; R_1, \ldots, R_{i-1}, R_i', R_{i+1}, \ldots, R_n) \not\ni z$.

(PR)  If $x, y \in C(S; \mathbb{R})$, $i \in N$, $R'_i \setminus \{x, y\}^2 = R_i \setminus \{x, y\}^2$,
and either $x\, E_i\, y$ and $x\, P'_i\, y$, or $y\, P_i\, x$ and $y\, E'_i\, x$,
then $y \notin C(S; R_1, \ldots, R_{i-1}, R'_i, R_{i+1}, \ldots, R_n)$.

Assume $x \in C(S; \mathbb{R})$, and now only one individual $i$ replaces either (i) a strict preference $y\, P_i\, x$ by an equivalence $y\, E'_i\, x$, or (ii) an equivalence $x\, E_i\, y$ by a strict preference $x\, P'_i\, y$, while all other individual preferences remain the same, giving a profile $(R_1, \ldots, R_{i-1}, R'_i, R_{i+1}, \ldots, R_n)$. Then the change should neither make $x$ unacceptable, nor make any former unacceptable alternative $z$ acceptable. This strong version (sNNR) of the well known condition of *non-negative responsiveness* is easily seen to hold for all of our majority/minority-based rules so far:

*Proof.* The relevant cardinalities change only for $x$ and $y$, and at most by one: either (i) $r_{xy}$ and $e_{xy}$ increase and $p_{yx}$ decreases, or (ii) $p_{xy}$ increases, and $e_{xy}$ and $r_{yx}$ decrease. Also, $d_S$ remains unchanged. Consequently, each of our majority/minority relations $A_\alpha$ is changed to some $A'_\alpha$ that differs only in that possibly the arrow $x \to y$ is added and/or the arrow $y \to x$ is removed. Therefore $x$ remains $\operatorname{tr}_S(A_\alpha)$-optimal. Now let $A_\alpha$ be one of those relations for which $z$ is not $\operatorname{tr}_S(A_\alpha)$-optimal, say $w \operatorname{tr}_S(A_\alpha) z$ but not $z \operatorname{tr}_S(A_\alpha) w$, and assume that $z$ becomes $\operatorname{tr}_S(A'_\alpha)$-optimal. (i) If still $w \operatorname{tr}_S(A'_\alpha) z$, the now also existing path $z\, A'_\alpha \cdots A'_\alpha\, w$ must contain the only possible new arrow $x \to y$, i.e. $z\, A'_\alpha \cdots A'_\alpha\, x\, A'_\alpha\, y\, A'_\alpha \cdots A'_\alpha\, w$. But then $z \operatorname{tr}_S(A_\alpha) x$, hence $w \operatorname{tr}_S(A_\alpha) x$, and thus $x \operatorname{tr}_S(A_\alpha) w$, since $x$ is $\operatorname{tr}_S(A_\alpha)$-optimal. This gives the contradiction $z \operatorname{tr}_S(A_\alpha) w$. (ii) If, on the other hand, no longer $w \operatorname{tr}_S(A'_\alpha) z$, the corresponding $A_\alpha$-path must have contained the now removed arrow $y \to x$, i.e. $w\, A_\alpha \cdots A_\alpha\, y\, A_\alpha\, x\, A_\alpha \cdots A_\alpha\, z$, thus also $x \operatorname{tr}_S(A'_\alpha) z$. But then $z \operatorname{tr}_S(A'_\alpha) x$, since $z$ is $\operatorname{tr}_S(A'_\alpha)$-optimal, hence also $z \operatorname{tr}_S(A_\alpha) x$. Together with $x \operatorname{tr}_S(A_\alpha) w$ (by optimality of $x$), we get again the contradiction $z \operatorname{tr}_S(A_\alpha) w$. □

The related condition of *positive responsiveness* (PR) expresses the idea that, if a "non-deterministic" choice $C(S)$ with $|C(S)| \geqslant 2$ is interpreted as a kind of "social undecidedness", it should be possible to turn it into a "deterministic" choice with $|C(S)| = 1$ by only slight changes in the individual preferences. However, the same objection that was given against (sP) is also a point against (PR), and none of our rules satisfies it, since for the three profiles

$$\begin{pmatrix} d & c & b \\ | & | & | \\ a & b & a \\ | & | & | \\ c & d & d \\ | & | & | \\ b & a & c \end{pmatrix}, \begin{pmatrix} d & c & \\ | & | & \\ a & b & ab \\ | & | & | \\ c & d & d \\ | & | & | \\ b & a & c \end{pmatrix}, \text{ and } \begin{pmatrix} d & c & a \\ | & | & | \\ a & b & b \\ | & | & | \\ c & d & d \\ | & | & | \\ b & a & c \end{pmatrix},$$

each of the majority relations has always the same set of optimals.

**Condorcet-type conditions.** It seems reasonable that if an alternative is strictly preferred to each other alternative by a (possibly varying) majority, it should be the one and only choice. Also, if $x$ is acceptable compared to every single alternative $y \in S$, one might want it to be acceptable in all of $S$, too (cf. [?]):

(C) $C(S) = \{x\}$ whenever $C(\{x,y\}) = \{x\}$ for all $y \in S$.

(GC) $x \in C(S)$ whenever $x \in C(\{x,y\})$ for all $y \in S$.

Until now, these are only fulfilled by some of our majority/minority-based rules. When only antisymmetric or only complete relations are used, (C) holds: Suppose we use the common optimal elements of a chain $\operatorname{tr}_S(F_1) \subseteq \cdots \subseteq \operatorname{tr}_S(F_m)$ of quasi-orders, and (i) either all of the relations $F_k$ are antisymmetric or (ii) all of them are complete. Then $\forall\, y \neq x\; \exists\, k : x\, F_k\, y\, \not\!\!F_k\, x$ implies $\forall\, y \neq x : x\, F_m\, y\, \not\!\!F_1\, x$, so that no $y \neq x$ is $\operatorname{tr}_S(F_m)$-optimal in case of (i) resp. $\operatorname{tr}_S(F_1)$-optimal in case of (ii). The variant (GC) is even fulfilled when both antisymmetric and complete (but no other) types are used: Suppose $x$ is $F_k|_{\{x,y\}}$-optimal for all $y \in S$. Then $y\, F_k\, x$ for no $y \in S \setminus \{x\}$ if $F_k$ is antisymmetric, and $x\, F_k\, y$ for all $y \in S$ if it is complete instead. In both cases $x$ is $\operatorname{tr}_S(F_k)$-optimal.

However, we can still also use rules based on types that are neither antisymmetric nor complete, since any rule $\mathcal{C}$ can be modified so as to fulfill (C) and/or (GC): Let $C^{(\mathrm{C})}(S)$ be (the uniquely determined) singleton $\{x\}$ if $C(\{x,y\}) := \{x\}$ for all $y \in S$, and $C^{(\mathrm{C})}(S) := C(S)$ otherwise. This cuts $\mathcal{C}$ down to its (C)-*modification* $\mathcal{C}^{(\mathrm{C})}$. Moreover,

$$C^{(\mathrm{GC})}(S) := C(S) \cup \{x \in S : x \in C(\{x,y\}) \text{ for all } y \in S\}$$

yields a (GC)-*modified* rule $\mathcal{C}^{(\mathrm{GC})}$ which gives larger sets than $\mathcal{C}$ in general.

**Lemma 4** *If $\mathcal{C}$ satisfies (wIm$_A$), (Im$_A$), (P), (sP), or (sP'), then so do $\mathcal{C}^{(\mathrm{C})}$ and $\mathcal{C}^{(\mathrm{GC})}$, respectively.*

*Proof.* As for (C): Under the assumption, $y\, A\, x \in C(\{x,y\})$ implies $x\, A\, y$. As for (GC): Analogously, $y\, A\, x \in C^{(\mathrm{C})}(S) \setminus C(S)$ implies $x \in C(\{x,y\})$, hence again $x\, A\, y$. □

Unfortunately, (sIm$_A$) is preserved by neither modification: (GC): For the profile

$$\begin{pmatrix} & & c & a \\ & & | & | \\ ac, & b, & b, & c \\ | & | & | & | \\ b & ac & a & b \end{pmatrix},$$

we get $M_{.5} = M_{.5}^S = B_{.5} = c \to b \leftrightarrow a$, so taking all $\operatorname{tr}_S(M_{.5})$-optimals leads to $a \in C(\{a,b\}) \cap C(\{a,c\})$, but $a$ is not optimal. (C): For

$$\begin{pmatrix} a & a & c & & & & \\ | & | & | & & & & \\ c, & c, & b, & b, & b, & ac, & ac \\ | & | & | & | & | & | & | \\ b & b & a & ac & ac & b & b \end{pmatrix},$$

we get $\{M_\alpha : \frac{1}{7} = \beta < \alpha \leqslant \gamma = \frac{4}{7}\} = \{c \leftrightarrow b \leftrightarrow a \to c,\ c \to b \leftrightarrow a,\ c \to b \leftarrow a\}$, so $\mathcal{M}_{(\beta,\gamma]}$ leads to $C(\{a,b\}) = C(\{a,c\}) = \{a\}$, hence $C^{(\mathrm{C})}(S) = \{a\}$, but $a$ is not $\operatorname{tr}_S(M_{3/7})$-optimal.

**Idempotency.** It would be somewhat strange if the same rule, applied to its first result, would "cut down" the choice set further instead of leaving it unchanged. One would rather expect that

(Id) *The choice function $S \mapsto C(S; \mathbb{R})$ is idempotent,*
*that is, $C(C(S; \mathbb{R}); \mathbb{R}) = C(S; \mathbb{R})$.*

Although, in the first place, of our rules only those based on complete relations satisfy (Id), any rule $\mathcal{C}$ has its (Id)-*modification* $\mathcal{C}^{(\mathrm{Id})}$, which is defined by $C^{(\mathrm{Id})}(S) := \bigcap_{m=1}^{\infty} C^m(S)$ and obviously fulfills (I), (Iso), (wIm$_A$), (Im$_A$), (sIm$_A$), (wP), (sP), (sP'), (NNP), (C), and (GC), if only $\mathcal{C}$ does.

**Ratio rules.** In [?], Pattanaik and Sengupta consider the special class SPCF of choice rules that are defined only on sets of two alternatives and for a certain kind of fuzzy relation instead of a preference profile. With the convention $0/0 := 1$, the quotients

$$\delta(x,y) := \frac{r_{xy}}{d_{xy}}, \quad \pi(x,y) := \frac{p_{xy}}{p_{xy} + p_{yx}} \quad \text{and} \quad \varrho(x,y) := \frac{r_{xy}}{r_{xy} + r_{yx}}$$

that correspond to $D_\alpha$, $P_\alpha$, and $R_\alpha$, respectively, are examples of such connected fuzzy relations, while the quotients that correspond to the other majority types are not connected. It turns out that the rules $\mathcal{B}_{(0,\gamma]}$, $\mathcal{M}_{(0,\gamma]}$, and $\mathcal{M}^S_{(0,\gamma]}$ with $\gamma \geqslant \frac{1}{2}$ give $x \in C(\{x,y\})$ if and only if $\pi(x,y) \geqslant \pi(y,x)$. In other words, the induced SPCFs are both *ratio* and *difference rules of type 1* in the sense of Pattanaik and Sengupta. The same is true for the rules $\mathcal{D}_{(\beta,1]}$, $\mathcal{N}_{(0,1]}$, and $\mathcal{N}^S_{(0,1]}$ with $\beta \leqslant \frac{1}{2}$ when we take $\varrho(x,y)$ instead of $\pi(x,y)$.

Moreover, $x$ is $\mathrm{tr}_{\{x,y\}}(P_\alpha)$-optimal if and only if $\pi(x,y) \geqslant \frac{\alpha}{1-\alpha}\pi(y,x)$ (for $\alpha \leqslant \frac{1}{2}$) resp. if $\pi(y,x) \leqslant c\pi(x,y)$ (for $\alpha > \frac{1}{2}$ and $c = \frac{\alpha}{1-\alpha} - \varepsilon \geqslant 1$ with some sufficiently small $\varepsilon > 0$). Therefore, also all rules $\mathcal{P}_{(\beta,\gamma]}$ are ratio rules, and the analogue holds for $\mathcal{R}_{(\beta,\gamma]}$ when using $\varrho(x,y)$ instead.

Furthermore, it might be interesting that, for a profile of total quasi-orders, the functions $d_P(x,y) := 1 - \pi(x,y)$ and $d_R(x,y) := 1 - \varrho(x,y)$ are distance functions (cf. [?]), that is, they fulfill the triangle inequality $d(x,y) + d(y,z) \geqslant d(x,z)$. In fact,

$$\frac{p_{xy}}{p_{xy} + p_{yx}} + \frac{p_{yz}}{p_{yz} + p_{zy}} - \frac{p_{xz}}{p_{xz} + p_{zx}} + 1$$
$$= \frac{p_{xy}}{p_{xy} + p_{yx}} + \frac{p_{yz}}{p_{yz} + p_{zy}} + \frac{p_{zx}}{p_{zx} + p_{xz}}$$
$$= \frac{p_{xy}}{n - e_{xy}} + \frac{p_{yz}}{n - e_{yz}} + \frac{p_{zx}}{n - e_{zx}}$$
$$\leqslant \frac{p_{xy} + e_{xy} - e_{xyz}}{n - e_{xyz}} + \frac{p_{yz} + e_{yz} - e_{xyz}}{n - e_{xyz}} + \frac{p_{zx} + e_{zx} - e_{xyz}}{n - e_{xyz}}$$
$$\leqslant \frac{(r_{xy} + r_{yz} + r_{zx}) - 3e_{xyz}}{n - e_{xyz}}$$
$$\leqslant \frac{(2(n - e_{xyz}) + 3e_{xyz}) - 3e_{xyz}}{n - e_{xyz}} = 2$$

and
$$\frac{r_{xy}}{r_{xy} + r_{yx}} + \frac{r_{yz}}{r_{yz} + r_{zy}} - \frac{r_{xz}}{r_{xz} + r_{zx}} + 1$$
$$= \frac{r_{xy}}{r_{xy} + r_{yx}} + \frac{r_{yz}}{r_{yz} + r_{zy}} + \frac{r_{zx}}{r_{zx} + r_{xz}}$$
$$= \frac{n - p_{yx}}{2n - (p_{yx} + p_{xy})} + \frac{n - p_{zy}}{2n - (p_{zy} + p_{yz})} + \frac{n - p_{xz}}{2n - (p_{xz} + p_{zx})}$$
$$\leqslant \frac{n + p_{xy} - \ell_{xyz}}{2n - \ell_{xyz}} + \frac{n + p_{yz} - \ell_{xyz}}{2n - \ell_{xyz}} + \frac{n + p_{zx} - \ell_{xyz}}{2n - \ell_{xyz}}$$
$$\leqslant \frac{3n + (n + \ell_{xyz}) - 3\ell_{xyz}}{2n - \ell_{xyz}} = 2,$$

where $e_{xyz}$ and $\ell_{xyz}$ are the numbers of individuals $i$ that have $x\,E_i\,y\,E_i\,z$ resp. a linear order on $x, y, z$. In fuzzy preference theory, this property of $\pi(x,y)$ and $\varrho(x,y)$ is also called $T_2$-transitivity (cf. [?,?,?]). For profiles of 3-acyclic relations (i.e., when $x\,P_i\,y\,P_i\,z\,P_i\,x$ for no $i, x, y, z$), obviously also $d_M(x,y) := 1 - \frac{p_{xy}}{n}$ fulfills the triangle inequality.

**Other conditions.**

(CA) $x \in C(S)$, if $x \in C(S')$, $x \in S \subseteq S'$, and $C(S) \cap C(S') \neq \emptyset$.

($\beta$) If $x, y \in C(S)$, and $S \subseteq S'$, then $x \in C(S') \iff y \in C(S')$.

(SUA) $C(S) = C(S')$ whenever $S \subseteq S'$ and $C(S') \subseteq C(S)$.

In [?], Sen gives the following intuition for a slightly stronger version of *Chernoff's condition* (CA)*:* "if $x$ is a best alternative in a given set and belongs to a certain subset of it, then $x$ must also be best in that subset, e.g. a world champion must also be champion in his country". Although, with a suitable interpretation of "best", this intuition is absolutely correct in my opinion, it is not, however, accurately expressed by (CA) or Sen's stronger condition ($\alpha$), since we cannot, in general, expect the elements of $C(S)$ to be "best" in a sense as strong as indicated by the word "champion". Consider the well-known minimal example for a cyclic majority relation,
$$\begin{pmatrix} a & b & c \\ | & | & | \\ b, & c, & a \\ | & | & | \\ c & a & b \end{pmatrix},$$
where $X := \{a, b, c\}$. Here, $C(X)$ should be $X$ indisputably, and, at least to me, it is almost equally obvious that $C(\{a,b\})$ should be $\{a\}$. But (CA) would require $C(\{a,b\}) = \{a,b\}$, which is absolute nonsense. However, when $x \in S \subseteq S'$ is a best alternative in $S'$ in the sense that $x \in C(\{x,y\})$ for all $y \in S'$, it should indeed be in $C(S)$ (and as well in $C(S')$, of course), but this requirement is already expressed in condition (GC).

As a support for ($\beta$), Sen gives the intuition that, "if two alternatives are both best in a certain subset, then one can be best for the whole set if and only if so is the other", which would be correct when "best" could is interpreted as "at least as good as any other". But, as we have already seen, the reason for $C(S)$ containing more than one alternative may be that some alternatives are just *not* comparable: Assume that $X = \{a, b, c, d\}$ and all individuals have the same preferences

$$\begin{array}{cc} a & c \\ | & | \\ b & d \end{array}.$$

Then (Iso), (I), and (P) require $C(X) = \{a, c\}$ and $C(\{a, d\}) = \{a, d\}$, which is a perfect solution, but not in accordance with ($\beta$).

A similar example shows that (SUA) is not absolutely reasonable in presence of undecidedness: Suppose both of two individuals have $x\,U_i\,y$, which results in $C(\{x, y\}) = \{x, y\}$ when (Iso) and (I) hold, and now a third alternative $z$ becomes feasible. If both are undecided about $z$ except that one prefers $x$ to $z$ and the other prefers $z$ to $y$, $x$ should become the only acceptable alternative, in contradiction to (SUA). Suzumura [**?**] shows that (CA) and (SUA) are implied by many other conditions which I consequently won't discuss here.

## 6 Invariance under decomposition into generalized components

When the number of alternatives is large, an important technique in actual decision processes is that the alternatives may be grouped together to form a smaller number of clusters, since then one can first choose between these clusters and then inside the chosen cluster(s). For example, if the alternatives are some thousand video taped movies in a store, it is practically impossible, even for a single person, to write down a preference relation by comparing each pair of movies separately. Instead, the alternatives are usually compared in classes, i.e. the group may classify them according to their genre, then choose some genre, and finally choose a single movie belonging to that genre. There is no problem in doing so, at least as long as all individuals consider a classification by genre also significant for their *own* preferences. But, in the above example, some individuals might prefer a clustering by age, length, or actors, so that they probably would not accept the suggested procedure to first choose a genre.

I will not discuss the difficult question of how, in general, to produce a meaningful clustering which all individuals will accept, but there is one situation I consider important for an axiomatic discussion of choice functions: Suppose there is a subset $B \subseteq X$ of alternatives that "belong together" in the sense that, for each individual $i$ and each alternative $x$ outside of $B$, $i$ values all alternatives in $B$ alike when compared with $x$. Then it seems very likely that all individuals would agree that $B$ is a meaningful cluster and can therefore be treated like a single alternative in the beginning. However, there might at the same time also exist other subsets $B'$ for which the very same is true. Therefore, it would be good to know that, in the end, the resulting choice did not depend on the use of classification $B$ or any other classification $B'$, but would have been the same without any classification at all.

To be able to make this argument more precise, let us extend the concept of a "component" of a tournament (cf. [?]) to the case of general binary relations: A *(generalized) component* of a reflexive relation $R$ on $X$ is a subset $B \subseteq X$ with $|B| \geqslant 2$ such that, for all $x \in X \setminus B$, whenever $x\,R\,y$ for *some* $y \in B$, $x\,R\,y$ for *all* $y \in B$, and whenever $y\,R\,x$ for some $y \in B$, $y\,R\,x$ for all $y \in B$. An equivalent formulation of the latter is this: For all $x \in X \setminus B$, $B \cap xR \neq \emptyset$ implies $B \subseteq xR$, and $B \cap Rx \neq \emptyset$ implies $B \subseteq Rx$, where I used the usual order-theoretic notation $xR = \{y \in X : x\,R\,y\}$ and $Rx = \{y \in X : y\,R\,x\}$. If $B$ is a common component of $R_1, \ldots, R_n$, then we shall call it a component of $\mathbb{R}$. The components of $R$ resp. $\mathbb{R}$ build a hull system, and if $B, B'$ are incomparable but intersecting components, also $B \cup B'$ and $B \setminus B'$ are components.[3]

To treat $B$ like a single alternative corresponds to using a new set of alternatives $S/B := S \setminus B \cup \{B\}$ and the *quotient profile* $\mathbb{R}/B := (R_1/B, \ldots, R_n/B)$ of individual preferences

$$R_i/B := R_i|_{X \setminus B} \cup \{(B,B)\}$$
$$\cup \big(\{B\} \times (BR_i \setminus B)\big) \cup \big((R_iB \setminus B) \times \{B\}\big),$$
$$\text{where } BR_i = \{y \in X : x\,R_i\,y \text{ for some } x \in B\}$$
$$\text{and } R_iB = \{y \in X : y\,R_i\,x \text{ for some } x \in B\}.$$

In other words, for all $x,y \in X/B$, $x\,R_i/B\,y$ if and only if either (i) $x \neq B$, $y \neq B$, and $x\,R_i\,y$, or (ii) $x = y = B$, or (iii) $x = B$, $y \neq B$, and $z\,R_i\,y$ for all $z \in B$, or (iv) $x \neq B$, $y = B$, and $x\,R_i\,z$ for all $z \in B$. For example:

$$\mathbb{R} = \begin{pmatrix} \begin{array}{cc} a & c \\ /\backslash & /\backslash \\ b\ c & bb'\ d \\ |\ | & |\ | \\ b'\ c' & c'\ a \\ | \\ d \end{array} \end{pmatrix}, B = \{b, b'\} \Longrightarrow \mathbb{R}/B = \begin{pmatrix} \begin{array}{cc} a & c \\ /\backslash & /\backslash \\ B\ c & B\ d \\ |\ | & |\ | \\ c' & c'\ a \\ | \\ d \end{array} \end{pmatrix}.$$

The set $\{c, c'\}$ is only a component of the first relation.

Now, our above consideration finally yields the following condition on $\mathcal{C}$:

(CC) *For each component $B \subseteq S$ of $\mathbb{R}$:*

$$C(S; \mathbb{R}) = \begin{cases} C(S/B; \mathbb{R}/B) & \text{if } B \notin C(S/B; \mathbb{R}/B) \\ C(S/B; \mathbb{R}/B) \setminus \{B\} \\ \quad \cup\, C(B; \mathbb{R}) & \text{if } B \in C(S/B; \mathbb{R}/B). \end{cases}$$

---

[3] One has to be a bit careful about the terminology here: In case of graphs or quasi-ordered sets, the generalized components are usually called "(lexicographic) factors" or "blocks" [?], respectively. In graph and order theory, on the other hand, (ordinal or cardinal) "components" are special cases of generalized components, while a "block" of a graph is yet another thing.

This expresses the claim that, whatever component $B$ we use, the acceptable alternatives are (i) those *outside* of $B$ that are acceptable when treating $B$ like a single alternative, and (ii) if $B$ itself is then acceptable, too, also those alternatives *inside* of $B$ that would be acceptable when only the $B$-alternatives were considered.

In [**?**], Laffond et al. introduced (CC) for the more special case of tournaments, but using partitions into components instead of a single component.

Let us define a partial order (*"finer than"*) between rules by setting $\mathcal{C} \leqslant \mathcal{C}'$ if $C(S;\mathbb{R}) \subseteq C'(S;\mathbb{R})$ for all $S, \mathbb{R}$. It is easily checked that whenever $(\mathcal{C}_i)_i$ is a family of rules with (CC) such that $C'(S;\mathbb{R}) := \bigcap_i C_i(S;\mathbb{R}) \neq \emptyset$ for all $\mathbb{R}$ and $S \neq \emptyset$, the rule $\mathcal{C}' : (S;\mathbb{R}) \mapsto C'(S;\mathbb{R})$ also fulfills (CC). Therefore, any rule $\mathcal{C}$ trivially possesses a finest rule $\mathcal{C}' \geqslant \mathcal{C}$ with (CC), called the (CC)-hull of $\mathcal{C}$. This $\mathcal{C}'$ might give considerably larger choice sets than $\mathcal{C}$, in many cases it will even give constantly $C'(S) = S$ (cf. [**?**]). Also, no simple algorithm to determine $C'(S)$ is known, and as there can be exponentially many components of $S$, it seems very likely that there can be no such algorithm with polynomial time complexity. Therefore, let us look at another way to make a rule fulfill (CC):

**Theorem 5** *Let $\mathcal{C} : (S, \mathbb{R}) \mapsto C(S;\mathbb{R})$ be a social choice rule that fulfills* (Iso)*, *(I)*, (GC)*, *and* (C)*. Then $C^{(\mathrm{CC})}(S;\mathbb{R}) :=$*

$$\begin{cases} C(S;\mathbb{R}) & \text{if no } B \subset S \text{ is a component of } \mathbb{R} \\ C^{(\mathrm{CC})}(S/B;\mathbb{R}/B) & \text{if } B \subset S \text{ is a component of } \mathbb{R} \\ & \text{ and } B \notin C^{(\mathrm{CC})}(S/B;\mathbb{R}/B) \\ C^{(\mathrm{CC})}(S/B;\mathbb{R}/B) \setminus \{B\} & \text{if } B \subset S \text{ is a component of } \mathbb{R} \\ \cup\, C^{(\mathrm{CC})}(B;\mathbb{R}) & \text{ and } B \in C^{(\mathrm{CC})}(S/B;\mathbb{R}/B), \end{cases}$$

*where $\subset$ denotes proper containment, recursively defines another rule $\mathcal{C}^{(\mathrm{CC})}$ such that:*

1. *$C^{(\mathrm{CC})}(S;\mathbb{R}) = C(S;\mathbb{R})$ if no component of $\mathbb{R}$ is properly contained in $S$.*
2. *$\mathcal{C}^{(\mathrm{CC})}$ fulfills* (Iso)*, (I), (GC), (C), and (CC).*

The somewhat lengthy proof is to be found in the appendix.

**Corollary 6** *When $C(S;\mathbb{R})$ can be computed in polynomial time for all $S$ and $\mathbb{R}$ then so can $C^{(\mathrm{CC})}(S;\mathbb{R})$.*

*Proof.* For $A \subseteq S$, $|S| = k$, and a profile $\mathbb{R}$ on $S$, the set

$$\varphi(A) := A \cup S \setminus \bigcap_{i \in N} \left( \bigcap_{x \in A} xP_i \cup \bigcap_{x \in A} xE_i \cup \bigcap_{x \in A} P_i x \cup \bigcap_{x \in A} xU_i \right)$$

can be computed in $O(k^2 n)$ time. A straightforward proof shows that, for $x, y \in S$, $x \neq y$, $B_{xy} := \varphi^k(\{x, y\})$ is the smallest component of $\mathbb{R}$ that contains $x$ and $y$; it can thus be computed in $O(k^3 n)$ time. If $\mathbb{R}$ has a proper component, one of the $\binom{k}{2}$ many components $B_{xy}$ must be proper and minimal, hence a minimal proper component $B$ can be found in $O(k^5 n)$ time. When taking always minimal components in the recursion, $C^{(\mathrm{CC})}(S)$ can be determined by computing $C(S_j)$ for at most $k-1$ different sets $S_j \subseteq S$, because all of the occurring components and the final quotient are isomorphic to subsets of $S$. Thus, the time complexity of $\mathcal{C}^{(CC)}$ is at most $O(k^6 n)$ times that of $\mathcal{C}$. □

We have already seen that the requirements of the theorem are fulfilled by all rules based on only antisymmetric or only complete majority/minority relations, and also by the (C)- and (GC)-modified rules.

**Lemma 7** *If components of $\mathbb{R}$ are always components of $A$, the modification $\mathcal{C} \mapsto \mathcal{C}^{(CC)}$ preserves* (wIm$_A$), (Im$_A$), *and* (sIm$_A$).

*Proof.* Assume that $B \subseteq S$ is a component of $\mathbb{R}$. Then $B$ is also a component of $A$. For (sIm$_A$), let $x \in S \setminus B$, $x \in C(S/B; \mathbb{R}/B)$, and $y \operatorname{tr}_S(A) x$. Then

(i) $y \notin B$: $y \operatorname{tr}_{S/B}(A/B) x \implies x \operatorname{tr}_{S/B}(A/B) y \implies x \operatorname{tr}_S(A) y$,
(ii) $y \in B$: $B \operatorname{tr}_{S/B}(A/B) x \implies x \operatorname{tr}_{S/B}(A/B) B \implies x \operatorname{tr}_S(A) y$.

On the other hand, let $B \in C(S/B; \mathbb{R}/B)$, $x \in C(B; \mathbb{R})$, and $y \operatorname{tr}_S(A) x$. Then

(iii) $y \notin B$: $y \operatorname{tr}_{S/B}(A/B) B \implies B \operatorname{tr}_{S/B}(A/B) y \implies x \operatorname{tr}_S(A) y$,
(iv) $y \in B$, $y \operatorname{tr}_B(A) x$: $\quad x \operatorname{tr}_B(A) y \implies x \operatorname{tr}_S(A) y$,

and (v) otherwise $y \in B$ and the $A$-path from $y$ to $x$ must leave $B$, i.e. $y \operatorname{tr}_S(A) y' A z \operatorname{tr}_S(A) w A x' \operatorname{tr}_S(A) x$ for some $x', y' \in B$ and $z, w \in S \setminus B$, which implies $x A z \operatorname{tr}_S(A) w A y$ and thus $x \operatorname{tr}_S(A) y$.

For (wIm$_A$) and (Im$_A$), the arguments are analogous, only that case (v) does not occur. □

**Lemma 8** *If $\mathcal{C}$ fulfills* (Iso) *and* (CC) *then so does* $\mathcal{C}^{(Id)}$.

*Proof.* Exactly as in [**?**], where the special case of tournaments was considered. □

## 7 Conclusion

As we have seen, it is not a real problem to allow individual preferences to be incomplete and even cyclic, i.e., any reflexive relation. In fact, Lemma 1 shows that we should at least allow all quasi-orders. The requirements of immunity from binary arguments, (Im$_A$), for a large class of majority and/or minority relations $A$ of different strength directly lead to social choice rules which already fulfill some of the most frequently discussed axioms. Using proper modifications, one may fulfill additional axioms: Lemmata 2, 3, 4, 7, 8 and Theorem 5 imply

**Corollary 9** *Let $F$ be a chain of majority relations with $M_1 \in F$, and $\mathcal{C}$ be the rule for which $C(S)$ is the set of common optimals of $\{\operatorname{tr}_S(A) : A \in F\}$. Then the rule $\mathcal{C}^{(C)(GC)(CC)(Id)}$ fulfills the axioms* (Iso), (I), (P), (C), (GC), (wIm$_A$), *and* (Im$_A$) *for all $A \in F$. Moreover, $N_1 \in F$ implies* (sP), *and $B_1 \in F$ implies* (sP').

Depending on whether (sP) or (sP') is required, a look at Figure 1 shows that one may for example choose one of the chains

$$F = \{D_\beta, \ldots, D_\gamma, N^S_\delta, \ldots, N^S_\varepsilon, N_\zeta, \ldots, N_\eta, N_1, M_1\}$$

with $0 < \beta \leqslant \gamma \leqslant \delta \leqslant \varepsilon \leqslant \zeta \leqslant \eta \leqslant 1$,

$$F = \{D_\beta, \ldots, D_\gamma, R_\delta, \ldots, R_\varepsilon, P_\zeta, \ldots, P_\eta, E_{\frac{1}{2}},$$
$$U_\vartheta, \ldots, U_\kappa, B_\lambda, \ldots, B_\mu, M^S_\nu, \ldots, M^S_\xi, M_\pi, \ldots, M_\varrho, M_1\}$$

with $0 < \beta \leqslant \cdots \leqslant \eta \leqslant \frac{1}{2} \leqslant \vartheta \leqslant \cdots \leqslant \varrho$, or

$$F = \{D_\beta, \ldots, D_\gamma, E_\delta, \ldots, E_\varepsilon, P_\zeta, \ldots, P_\eta,$$
$$R_\vartheta, \ldots, R_\kappa, B_\lambda, \ldots, B_\mu, M^S_\nu, \ldots, M^S_\xi, M_\pi, \ldots, M_\varrho, M_1\}$$

with $0 < \beta \leqslant \frac{n}{2n-1} < \gamma \leqslant \cdots \leqslant \varrho$.

Furthermore, it was argued that considerations about possible power of minorities and the need to "protect" minorities against majorities may help deciding what types of relations and what range of strengths actually to use, where the types $B$ and $D$ seem to have the greatest appeal.

Many things remain to be done: First of all, I did not check whether the (C)- and (GC)-modifications also preserve (sNNR). Secondly, what about strategic voting? And, what seems to be most important to me, there should be some "experimental" investigations by means of stochastic simulations which would randomly generate profiles of preferences of different type: either arbitrary, or complete and/or acyclic and/or transitive and/or antisymmetric. By their means it should be possible to distinguish between the rules where the axiomatic approach alone does not provide a satisfactory evaluation.

In fact, first "on-the-fly"-simulations seem to suggest that, depending on the degrees of antisymmetry and completeness in the individual preferences and on the relative sizes of $S$ and $N$, either the combination of strict types like $B$ and $P$ or the combination of non-strict types like $D$ and $R$ has a very large probability of giving a "deterministic" choice $C(S)$ with $|C(S)| = 1$. Particularly, for $A \in \{R, P, E, D\}$ and $n \to \infty$, the probability that $|\{A_\alpha : \frac{1}{2} < \alpha \leqslant 1\}| = \binom{|S|}{2}$ tends to one, and in this case the rules $\mathcal{A}_{(.5,1]}$ must choose deterministically: Otherwise, assuming $x, y \in C(S)$, $x \neq y$, and, without loss of generality, $x A_{.5+\varepsilon} y$ because of (ii), one would conclude that also $y \operatorname{tr}_S(A_{.5+\varepsilon}) x$. Then (ii) implies that there is some $\alpha$ for which only one of $x \operatorname{tr}_S(A_\alpha) y$ and $y \operatorname{tr}_S(A_\alpha) x$ is true — a contradiction. That is, quite contrary to the simple majority rule, these rules lead to deterministic choices with a probability that tends to one for large $n$. This might also indicate that type $D$ surpasses $B$ here.

In case of large $S$, on the other hand, it is, unfortunately, a nontrivial task to generate a uniformly distributed random sample of acyclic or transitive relations on $S$ (see, for example, [?], and [?] for the case of order relations). For a simulation, one would therefore first have to search for a probability distribution on the possible preference relations that models actual distributions sufficiently well and at the same time allows for a quick generation of samples.

Finally, how can an actual implementation of these rules look like? While for only three or four alternatives, a group can perform the necessary steps by hand, this becomes quickly unworkable with more alternatives. That is, it might be quite difficult to design a *procedural, intuitive, plausible,* and *manageable* algorithm in the sense of Haake, Raith, and Su [**?**,**?**] (see also [**?**]). Probably the application of any sophisticated choice rule that allows for other individual preference relations than only total (quasi-)orders will have to involve a human or artificial moderator who first determines the individual preferences either by pairwise ballots or by interviewing each individual separately, then performs the computation of $C(S)$, and finally explains and defends this suggested result. This motivates the search for an interviewing technique that can determine an arbitrary individual preference relation on $S$ in a considerably smaller expected time than $\binom{|S|}{2}$, the time required for the pairwise ballots.

**Appendix**

*Proof of Theorem 5.* Mainly, we have to verify that $\mathcal{C}' := \mathcal{C}^{(\mathrm{CC})}$ is indeed well-defined, i.e., that the definition is not circular and does not depend on the choice of $B$. We can do this by induction on the size of $S$: Let $\mathbb{R}$ and $S$ be given. For $|S| = 1$, no component is properly contained in $S$, and the five conditions in 2. hold, because here $C'(S) = S$, independently of any profile. Now assume that well-definedness of $C'(\tilde{S}; \tilde{\mathbb{R}})$ and the five conditions have already been proved for all $\tilde{\mathbb{R}}$ and all $|\tilde{S}| < |S|$. Let $B_1, B_2 \subset S$ be components of $\mathbb{R}$. Then $|B_1|, |B_2|, |S/B_1|, |S/B_2| < |S|$, so that the definition is not circular, and, without loss of generality, we have one of the following five situations:

(i) $B_1 \cap B_2 = \emptyset$: Then $B_1$ and $B_2$ are also components of $\mathbb{R}/B_2$ and $\mathbb{R}/B_1$, respectively, the profile $\mathbb{R}_{21} := (\mathbb{R}/B_2)/B_1$ on $S_{21} := (S/B_2)/B_1$ is isomorphic to the profile $\mathbb{R}_{12} := (\mathbb{R}/B_1)/B_2$ on $S_{12} := (S/B_1)/B_2$, $\mathbb{R}/B_1|_{B_2} \cong \mathbb{R}|_{B_2}$, and $\mathbb{R}/B_2|_{B_1} \cong \mathbb{R}|_{B_1}$. Therefore, by induction, $C'(S; \mathbb{R}) \setminus (B_1 \cup B_2) = C'(S_{21}; \mathbb{R}_{21}) = C'(S_{12}; \mathbb{R}_{12})$. Moreover, when $B_1$ is used in the recursion, either

$$C'(S; \mathbb{R}) \cap B_1 = C'(B_1; \mathbb{R})$$

or it is empty, the former if $B_1 \in C'(S/B_1; \mathbb{R}/B_1)$. If, on the other hand, $B_2$ is used, either

$$C'(S; \mathbb{R}) \cap B_1 = C'(S/B_2; \mathbb{R}/B_2) \cap B_1 = C'(B_1; \mathbb{R}/B_2)$$

or it is empty, the former if $B_1 \in C'(S_{21}; \mathbb{R}_{21})$. But, by induction,

$$B_1 \in C'(S/B_1; \mathbb{R}/B_1) \iff B_1 \in C'(S_{12}; \mathbb{R}_{12}) \iff B_1 \in C'(S_{21}; \mathbb{R}_{21}),$$

and $C'(B_1; \mathbb{R}) = C'(B_1; \mathbb{R}/B_2)$. The case of $C'(S; \mathbb{R}) \cap B_2$ is analogous. This shows that $C'(S; \mathbb{R})$ does not depend on whether $B_1$ or $B_2$ is used in the recursion.

(ii) $B_1 \subseteq B_2$: Then $B_1$ is a component of $\mathbb{R}|_{B_2}$, $B_2/B_1$ is a component of $\mathbb{R}/B_1$, and $\mathbb{R}/B_2 \cong (\mathbb{R}/B_1)/(B_2/B_1)$. Again,

$$C'(S; \mathbb{R}) \setminus B_2 = C'(S/B_2; \mathbb{R}/B_2) = C'(S_{12}; \mathbb{R}_{12}).$$

Now, when directly using $B_2$,

$$C'(S;\mathbb{R}) \cap B_1 = C'(B_2;\mathbb{R}) \cap B_1 = C'(B_1;\mathbb{R})$$

if $B_2 \in C'(S/B_2;\mathbb{R}/B_2)$ and $B_1 \in C'(B_2/B_1;\mathbb{R}/B_1)$. On the other hand, when first $B_1$ is used,

$$C'(S;\mathbb{R}) \cap B_1 = C'(B_1;\mathbb{R})$$

if $B_1 \in C'(S/B_1;\mathbb{R}/B_1)$. But $B_1 \in B_2/B_1$ which is a component of $\mathbb{R}/B_1$, hence $B_1 \in C'(S/B_1;\mathbb{R}/B_1)$ is equivalent to $B_2 \in C'(S/B_2;\mathbb{R}/B_2)$ and $B_1 \in C'(B_2/B_1;\mathbb{R}/B_1)$.

Moreover, when we directly use $B_2$,

$$C'(S;\mathbb{R}) \cap (B_2 \setminus B_1) = C'(B_2;\mathbb{R}) \setminus B_1 = C'(B_2/B_1;\mathbb{R}/B_1) \setminus \{B_1\}$$

if $B_2 \in C'(S/B_2;\mathbb{R}/B_2)$. When $B_1$ is used instead,

$$C'(S;\mathbb{R}) \cap (B_2 \setminus B_1) = C'(S/B_1;\mathbb{R}/B_1) \cap (B_2 \setminus \{B_1\})$$
$$= C'(B_2/B_1;\mathbb{R}/B_1) \setminus \{B_1\}$$

if $B_2/B_1 \in C'((S/B_1)/(B_2/B_1);(\mathbb{R}/B_1)/(B_2/B_1))$. Again, this last condition is equivalent to $B_2 \in C'(S/B_2;\mathbb{R}/B_2)$, because there is a canonical isomorphism between $(\mathbb{R}/B_1)/(B_2/B_1)$ and $\mathbb{R}/B_2$ that maps $B_2/B_1$ to $B_2$.

(iii) $|B_1 \cap B_2| \geqslant 2$: Then also $B_3 := B_1 \cap B_2$ is a component of $\mathbb{R}$, so that we can apply (ii) two times because of $B_3 \subseteq B_1$ and $B_3 \subseteq B_2$.

(iv) $|B_1 \setminus B_2|, |B_2 \setminus B_1| \geqslant 2$: The reader may easily verify that then also $B_1' := B_1 \setminus B_2$ and $B_2' := B_2 \setminus B_1$ are components of $\mathbb{R}$, so that we may first apply (ii) on $B_1$ and $B_1'$, then (i) on $B_1'$ and $B_2'$, and finally (ii) on $B_2'$ and $B_2$.

(v) $|B_1| = 2$, $|B_1 \cap B_2| = 1$, say $B_1 = \{x,y\}$ and $x \in B_2$: As for $x$: Using $B_2$,

$$x \in C'(S;\mathbb{R}) \iff x \in C'(B_2;\mathbb{R}) \text{ and } B_2 \in C'(S/B_2;\mathbb{R}/B_2),$$

while using $B_1$,

$$x \in C'(S;\mathbb{R}) \iff x \in C'(B_1;\mathbb{R}) \text{ and } B_1 \in C'(S/B_1;\mathbb{R}/B_1).$$

But (I) and (Iso) show that both conditions are equivalent, because the canonical isomorphism between $\mathbb{R}|_{B_2}$ and $\mathbb{R}/B_1|_{S/B_1}$ maps $x$ to $B_1$, and the one between $\mathbb{R}|_{B_1}$ and $\mathbb{R}/B_2|_{S/B_2}$ maps $x$ to $B_2$.

As for $y$: When we use $B_2$,

$$y \in C'(S;\mathbb{R}) \iff y \in C'(S/B_2;\mathbb{R}/B_2).$$

When $B_1$ is used instead,

$$y \in C'(S;\mathbb{R}) \iff y \in C'(B_1;\mathbb{R}) \text{ and } B_1 \in C'(S/B_1;\mathbb{R}/B_1).$$

Again, both conditions are equivalent: First, $\mathbb{R}/B_2|_{S/B_2} \cong \mathbb{R}|_{B_1}$, so that (I) and (Iso) imply the equivalence of $y \in C'(S/B_2;\mathbb{R}/B_2)$ and $y \in C'(B_1;\mathbb{R})$. Moreover, because $B_2$ and $B_1$ are components, the latter implies that, for all $z \in B_2$, $y \in C'(\{y,z\}; \mathbb{R})$ and therefore $x \in C'(\{x,z\}; \mathbb{R})$. Now (GC) shows that then also $x \in C'(B_2; \mathbb{R})$, i.e. $B_1 \in C'(S/B_1; \mathbb{R}/B_1)$.

As for $z \in B_2 \setminus \{x\}$: Using $B_2$,

$$z \in C'(S; \mathbb{R}) \iff z \in C'(B_2; \mathbb{R}) \text{ and } B_2 \in C'(S/B_2; \mathbb{R}/B_2),$$

while if we use $B_1$,

$$z \in C'(S; \mathbb{R}) \iff z \in C'(S/B_1; \mathbb{R}/B_1).$$

As above, $\mathbb{R}/B_1|_{S/B_1} \cong \mathbb{R}|_{B_2}$, so that (I) and (Iso) imply the equivalence of $z \in C'(S/B_1; \mathbb{R}/B_1)$ and $z \in C'(B_2; \mathbb{R})$. Moreover, the latter implies that $C'(B_2; \mathbb{R}) \neq \{x\}$, so that, by (C), there must be some $w \in B_2 \setminus \{x\}$ such that $w \in C'(\{x,w\}; \mathbb{R})$. Since $B_1$ and $B_2$ are components, $w \in C'(\{y,w\}; \mathbb{R})$ and thus $x \in C'(\{y,x\}; \mathbb{R})$, so that finally $B_2 \in C'(S/B_2; \mathbb{R}/B_2)$.

Now a straightforward proof shows that (C) and (GC) also hold for $S$ and $\mathbb{R}$, while (Iso) and (I) are obvious, and (CC) holds by definition. □